\documentstyle[12pt]{article}

\input amssym.def
\input amssym

\font\tenmsb=msbm10 scaled \magstep1
\font\sevenmsb=msbm7 scaled \magstep1
\font\fivemsb=msbm5 scaled \magstep1
\newfam\msbfam
\textfont\msbfam=\tenmsb
\scriptfont\msbfam=\sevenmsb
\scriptscriptfont\msbfam=\fivemsb
\def\Bbb#1{{\fam\msbfam\relax#1}}

\font\teneufm=eufm10 scaled \magstep1
\font\seveneufm=eufm7 scaled \magstep1
\font\fiveeufm=eufm5 scaled \magstep1
\newfam\eufmfam
\textfont\eufmfam=\teneufm
\scriptfont\eufmfam=\seveneufm
\scriptscriptfont\eufmfam=\fiveeufm
\def\frak#1{{\fam\eufmfam\relax#1}}

\def\DEMONSTRATION{\smallskip{\noindent{\it D{\'e}monstration}.\ }}

\newtheorem{proposition}{{\sc Proposition}}[section]

\newtheorem{corollaire}[proposition]{{\sc Corollaire}}

\newtheorem{lemme}[proposition]{{\sc Lemme}}
\newtheorem{theoreme}{{\sc Th{\'e}or{\`e}me}}

\def\carre{\square\bigskip}

\let\[\lbrack
\let\]\rbrack
\let\<\langle
\let\>\rangle
\let\rta\rightarrow

\def\NN{{\Bbb N}}
\def\ZZ{{\Bbb Z}}

\def\FF{{\Bbb F}}

\def\RR{{\rm R}}
\def\Ga{{{\Bbb G}_a}}

\def\Ql{{\bar{\Bbb Q}}_\ell}
\def\A{{\cal A}}
\def\L{{\cal L}}

\def\OO{{\cal O}}

\def\R{{\cal R}}

\def\Id{{\rm Id}}

\def\N{{\rm N}}

\def\pr{{\rm pr}}

\def\res{{\rm res\,}}

\def\d{{\rm d}}
\def\diag{{\rm diag\,}}
\def\GL{{\rm GL}}

\def\vp{\varpi}

\def\hfl#1#2{\smash{\mathop{\hbox to 12mm{\rightarrowfill}}
\limits^{\scriptstyle#1}_{\scriptstyle#2}}}
\def\vfl#1#2{\llap{$\scriptstyle #1$}\left\downarrow
\vbox to 6mm{}\right.\rlap{$\scriptstyle #2$}}
\def\hgfl#1#2{\smash{\mathop{\hbox to 12mm{\leftarrowfill}}
\limits^{\scriptstyle#1}_{\scriptstyle#2}}}

\def\diagramme#1{\def\normalbaselines{\baselineskip=1pt
\lineskip=6pt\lineskiplimit=1pt} \matrix{#1}}

\begin{document}
\title{Preuve d'une conjecture de Frenkel-Gaitsgory-Kazhdan-Vilonen}
\author{Ng{\^o} Bao Ch{\^a}u}
\date{}
\maketitle

\begin{abstract}
We prove a conjecture of Frenkel-Gaitsgory-Kazhdan-Vilonen
on some exponential sums related to the
geometric Langlands correspondence. Our main ingredients are the
resolution of Lusztig scheme of lattices introduced by Laumon 
and the decomposition theorem of Beilinson-Bernstein-Deligne-Gabber. 
\end{abstract}

\section{L'{\'e}nonc{\'e}}

Soient $k=\FF_q$ un corps fini, ${\cal O}=k[[\vp]]$ le corps des s{\'e}ries
formelles {\`a} une variable $\vp$ et $F$ son corps des fractions.
Soient $d$ et $n$ deux entiers naturels.
A la suite de Lusztig (\cite{Lus}), consid{\'e}rons le sch{\'e}ma $X_d$
de type fini sur $k$ dont l'ensemble
des $k$-points est celui des r{\'e}seaux $\R\subset\OO^n$ tels que
$\dim(\OO^n/\R)=d$. L'action de $\GL(n,\OO)$ sur l'ensemble de ces r{\'e}seaux
peut {\^e}tre vue comme l'action d'un groupe alg{\'e}brique $G_d$ avec
$G_d(k)=\GL(n,\OO/\vp^d\OO)$, sur $X_d$.

Les orbites de cette action
sont en nombre fini. Pour chaque $n$-partition
$\lambda$ de $d$, c'est-{\`a}-dire $\lambda=(\lambda_1\geq\cdots\geq\lambda_n\geq 0)$ 
avec $|\lambda|=\lambda_1+\cdots+\lambda_n=d$, notons $X_\lambda$ l'orbite de $G_d$
passant par le r{\'e}seau $\vp^\lambda\OO^n$ o{\`u}
$\vp^\lambda$ d{\'e}signe la matrice diagonale
$\diag(\vp^{\lambda_1},\ldots,\vp^{\lambda_n})$.
On a la stratification en parties localement ferm{\'e}es
$X=\bigcup_{|\lambda|=n} X_\lambda$ qui refl{\`e}te la
d{\'e}composition de Cartan
$$\GL(n,F)=\coprod_{\lambda_1\geq\cdots\geq\lambda_n}
\GL(n,\OO)\vp^\lambda \GL(n,\OO).$$
En effet, on a
$$X_\lambda(k)=\GL(n,\OO)\vp^\lambda\GL(n,\OO)/\GL(n,\OO).$$

Pour chaque $\lambda$, notons ${\bar X}_\lambda$ l'adh{\'e}rence de l'orbite
$X_\lambda$ dans $X_d$. Rappelons que $X_\mu\subset{\bar X}_\lambda$
si et seulement si $\mu\leq\lambda$ selon l'ordre partiel habituel
entre les $n$-partitions de $d$ :
$$\mu_1+\cdots+\mu_i\leq\lambda_1+\cdots+\lambda_i$$
pour tout $i=1,\ldots,n-1$ (\cite{Lus}).

Fixons un nombre premier $\ell$ diff{\'e}rent de la caract{\'e}ristique $p$ de $k$.
Soit $\Ql$ une cl{\^o}ture alg{\'e}brique de ${\Bbb Q}_\ell$.
Notons $\A_\lambda$ le complexe d'intersection $\ell$-adique
de ${\bar X}_\lambda$.

Pour chaque $\alpha\in\NN^{\,n}$ tel que $|\alpha|=d$, notons $S_\alpha$ la
partie localement ferm{\'e}e de $X_d$ dont l'ensemble des $k$-points
est celui des r{\'e}seaux $\R\subset\OO^n$ tels que pour tout $i$, on a
$$(\R\cap\bigoplus_{j=1}^i e_j\OO)/(\R\cap\bigoplus_{j=1}^{i-1} e_j\OO)
=(\bigoplus_{j=1}^{i-1}e_j\OO\oplus\vp^{\alpha_i}e_i\OO)/
\bigoplus_{j=1}^{i-1}e_j\OO$$
o{\`u} $(e_i)$ d{\'e}signe la base standard de $\OO^n$.
La stratification $X_d=\bigcup_{|\alpha|=d}S_\alpha$ refl{\`e}te la
d{\'e}composition d'Iwasawa
$$\GL(n,F)=\coprod_{\alpha\in\ZZ^{\,n}}N(F)\vp^\alpha\GL(n,\OO),$$
o{\`u} $N$ d{\'e}signe le sous-groupe des matrices triangulaires sup{\'e}rieures
unipotentes de $\GL(n)$. En effet, on a
$$S_\alpha(k)=N(F)\vp^\alpha\GL(n,\OO)/\GL(n,\OO).$$

La fonction trace de Frobenius de $\A_\lambda$ s'identifie naturellement
{\`a} une fonction $A_\lambda$ {\`a} support compact dans $\GL(n,F)$
qui est bi-$\GL(n,\OO)$-invariante.
Fixons un caract{\`e}re additif non trivial $\psi:k\rta\Ql^\times$
et notons $\theta:N(F)\rta\Ql^\times$ le caract{\`e}re d{\'e}fini par
$$\theta(n)=\psi(\sum_{i=1}^{n-1}\res(n_{i,i+1}\d\vp)).$$
Consid{\'e}rons l'int{\'e}grale
$$I(\vp^\alpha,A_\lambda)=\int_{N(F)}A_\lambda(n\vp^\alpha)\theta(n)\d n,$$
o{\`u} la mesure de Haar normalis{\'e}e $\d n$ de $N(F)$ attribue {\`a}
$N(\OO)$ la mesure $1$.
Dans \cite{FGKV}, Frenkel, Gaitsgory, Kazhdan et Vilonen ont d{\'e}montr{\'e} le
th{\'e}or{\`e}me suivant.

\begin{theoreme}
Si $\alpha\not=\lambda$, on a
$$I(\vp^\alpha,A_\lambda)=0.$$
Si $\alpha=\lambda$, on a
$$I(\vp^\lambda,A_\lambda)=q^{\<\lambda,\delta \>}$$
o{\`u}
$$\delta={1\over 2}(n-1,n-3,\ldots,1-n)$$
et o{\`u}
$$\<\lambda,\delta\>=\sum_{i=1}^n\lambda_i\delta_i.$$
\end{theoreme}  

Lorsque la suite $\alpha=(\alpha_1,\cdots,\alpha_n)$ n'est pas
d{\'e}croissante,
on peut trouver $n'\in N(F)\cap \vp^\alpha\GL(n,\OO)\vp^{-\alpha}$
tel que $\theta(n')\not=1$.
Or, comme $A_\lambda$ est bi-$\GL(n,\OO)$-invariante, on a
\begin{eqnarray*}
\int_{N(F)} A_\lambda(n\vp^\alpha)\theta(n)\d n & =&
\int_{N(F)} A_\lambda(nn'\vp^\alpha)\theta(n)\d n \cr
& = &\theta(n')^{-1}\int_{N(F)}A_\lambda(n\vp^\alpha)\theta(n)\d n \cr
\end{eqnarray*}
donc $I(\vp^\alpha,A_\lambda)=0$.
Le cas int{\'e}ressant est donc celui o{\`u} $\alpha$ est une $n$-partition de
$d$. Dans ce cas, on a
$$N(F)\cap\vp^{\alpha}\GL(n,\OO)\vp^{-\alpha}\subset N(\OO)$$
si bien que le caract{\`e}re $N(F)\rta k$ d{\'e}fini par
$$n\mapsto\sum_{i=1}^{n-1}\res(n_{i,i+1}\d\vp)$$
induit un morphisme $h_\alpha:S_\alpha\rta\Ga$.
Frenkel, Gaitsgory, Kazhdan et Vilonen ont conjectur{\'e} dans \cite{FGKV}
l'{\'e}nonc{\'e} suivant.

\begin{theoreme}
Si $\alpha\not=\lambda$, on a
$$\RR\Gamma_c(S_\alpha\otimes_k{\bar k},
\A_\lambda\otimes h_\alpha^*\L_\psi)=0.$$
Si $\alpha=\lambda$, on a
$$\RR\Gamma_c(S_\alpha\otimes_k{\bar k},
\A_\lambda\otimes h_\alpha^*\L_\psi)
=\Ql[-2\<\lambda,\delta \>](-\<\lambda,\delta \>).$$
Ici, ${\bar k}$ d{\'e}signe une cl{\^o}ture alg{\'e}brique de $k$ et $\L_\psi$
le faisceau d'Artin-Schreier sur ${\Bbb G}_{a,k}$ associ{\'e} {\`a} $\psi$.
\end{theoreme}

On peut d{\'e}duire de cet {\'e}nonc{\'e} g{\'e}om{\'e}trique le th{\'e}or{\`e}me de
Frenkel-Gaitsgory-Kazhdan-Vilonen cit{\'e} plus haut, via
la formule des traces de Gro\-thendieck.

Voici les grandes lignes de la d{\'e}monstration du th{\'e}or{\`e}me $2$.

On consid{\`e}re d'abord le cas plus facile $\alpha=\lambda$.
On d{\'e}montre que si $\mu<\lambda$
l'intersection $S_\lambda\cap X_\mu$ est vide
si bien que celle de $S_\lambda$ avec le support de $\A_\lambda$
est incluse dans $X_\lambda$. On d{\'e}montre aussi que
$S_\lambda\cap X_\lambda$ est un espace affine et que le morphisme
$h_\alpha$ restreint {\`a} $S_\lambda\cap X_\lambda$ est constant {\`a} valeur
$0$ d'o{\`u} le r{\'e}sultat dans le cas $\alpha=\lambda$.
C'est le contenu de la section 2.

Pour d{\'e}montrer l'assertion concernant le cas $\alpha\not=\lambda$,
on utilise la r{\'e}solution suivante du sch{\'e}ma $X_d$.
Cette r{\'e}solution a {\'e}t{\'e} introduite par Laumon dans un contexte
l{\'e}g{\`e}rement diff{\'e}rent (\cite{Lau}).
Soit ${\tilde X}_d$ le sch{\'e}ma de type fini sur $k$ dont l'ensemble
des $k$-points est celui des drapeaux de r{\'e}seaux
$$\OO^n=\R_0\supset\R_1\supset\cdots\supset\R_d=\R$$
tels que $\dim(R_{i-1}/\R_i)=1$.
Le morphisme $\pi:{\tilde X}_d\rta X_d$ d{\'e}fini par
$$(\R_0\supset\R_1\supset\cdots\supset\R_n)\mapsto\R_n$$
est une r{\'e}solution semi-petite
au sens de Goresky et MacPherson.
De plus, elle est {\'e}quivariante relativement {\`a} l'action de $G_d$
si bien qu'on a 
$$\RR\pi_*\Ql[\dim(X_d)]({1\over 2}\dim(X_d))
=\bigoplus_\lambda \A_\lambda\boxtimes V_\lambda$$
o{\`u} les $V_\lambda$ sont des $\Ql$-espaces vectoriels,
gr{\^a}ce au th{\'e}or{\`e}me de d{\'e}composition (\cite{BBD})
et {\`a} ce que les sous-groupes stabilisateurs dans $G_d$ sont tous
g{\'e}om{\'e}triquement connexes.

Par comparaison avec la construction de Lusztig de la correspondance de
Springer, on voit que $V_\lambda$ est l'espace de la repr{\'e}sentation du groupe
sym{\'e}trique ${\frak S}_d$ correspondant {\`a} la partition $\lambda$ de $d$
(\cite{Lus},\cite{BM}).
On utilisera seulement le fait que la dimension $V_\lambda$ est {\'e}gal au nombre
de $\lambda$-tableaux standards.

Il suffit clairement de d{\'e}montrer que
$$\displaylines{\RR\Gamma_c(S_\lambda\otimes_k{\bar k},\RR\pi_*\Ql\otimes h_\lambda^*\L_\psi)\cr
=V_\lambda[-2\<\lambda,\delta\>-d(n-1)]
(-\<\lambda,\delta\>-{1\over 2}d(n-1)).}$$
Pour cela, on {\'e}tudie la g{\'e}om{\'e}trie de
${\tilde S}_\lambda=S_\lambda\times_{X_d}{\tilde X}_d$.
On a
$$\RR\Gamma_c(S_\lambda\otimes_k{\bar k},\RR\pi_*\Ql\otimes h_\lambda^*\L_\psi)
=\RR\Gamma_c({\tilde S}_\lambda\otimes_k{\bar k},{\tilde h}_\lambda^*\L_\psi)$$
o{\`u} ${\tilde h}_\lambda$ est le morphisme compos{\'e}
$h_\lambda\circ(\pi|_{{\tilde S}_\lambda})$.

On d{\'e}montre que ${\tilde S}_\lambda$ est une r{\'e}union disjointe 
de parties localement ferm{\'e}es ${\tilde S}_\tau$ qui sont des espaces
affines de m{\^e}me dimension 
$$\<\lambda,\delta\>+{1\over 2}d(n-1)=\<\lambda,(n-1,\ldots,1,0)\>$$
o{\`u} $\tau$
parcourt l'ensemble des suites $(\alpha^i)_{i=0}^d$ avec
$\alpha^i=(\alpha^i_j)_{j=1}^n\in\NN^{\,n}$ v{\'e}rifiant
\begin{itemize}
\item $\alpha^{i-1}_j\leq\alpha^i_j$
pour $i=1,\ldots,d$ et pour $j=1,\ldots,n$ ;
\item $|\alpha^i|=\sum_{j=1}^n\alpha^i_j=i$ pour $i=0,\ldots,d$ ;
\item $\alpha^d=\lambda$. 
\end{itemize}

Si l'une de ces suites $\alpha^i$ n'est pas d{\'e}croissante, on d{\'e}montre
comme dans le cas {\'e}voqu{\'e} plus haut o{\`u} $\lambda$ n'est pas d{\'e}croissante,
que
$$\RR\Gamma_c(S_\tau\otimes_k{\bar k},{\tilde h}_\lambda^*\L_\psi)=0.$$
Les $\tau$ dont les membres $\alpha_i$ sont tous des suites d{\'e}croissantes d'entiers naturels,
correspondent bijectivement aux $\lambda$-tableaux standards.
C'est le contenu de la section 3.

\section{Etude de $S_\alpha$}

Pour tout $\alpha\in\NN^n$, $S_\alpha$ est isomorphe {\`a} un espace
affine dont on peut construire les coordonn{\'e}es explicites {\`a} l'aide
de l'uniformisante $\vp$. Notons ${\bar\OO}=\OO\otimes_k{\bar k}$
et ${\bar F}=F\otimes_k{\bar k}$.

\begin{lemme}
Pour tout r{\'e}seau $\R\in S_\alpha({\bar k})$, il existe une unique matrice
triangulaire sup{\'e}rieure de la forme
$$x=\pmatrix{\vp^{\alpha_1} & x_{1,2} & \cdots & x_{1,n}\cr
                & \vp^{\alpha_2} & \cdots & x_{2,n}\cr
                &  & \ddots &\vdots \cr
                & & & \vp^{\alpha_n}}$$
o{\`u} les $x_{i,j}$ sont des polyn{\^o}mes en $\vp$ {\`a} coefficients
dans ${\bar k}$ de degr{\'e} strictement inf{\'e}rieur
{\`a} $\alpha_i$, telle que $\R=x{\bar{\cal O}}^n$.
\end{lemme}

\DEMONSTRATION
Du fait que $\R\in S_\alpha({\bar k})$, il se d{\'e}compose en
$$\R=\R'\oplus (\vp^{\alpha_n}e_n+y){\bar{\cal O}}$$
o{\`u}
$$\R'=\R\cap\bigoplus_{j=1}^{n-1}e_j{\bar{\cal O}}\in S_{\alpha'}({\bar k})$$
avec $\alpha'=(\alpha_1,\ldots,\alpha_{n-1})$
et o{\`u} $y\in \bigoplus_{j=1}^{n-1}e_j{\bar{\cal O}}$ est bien d{\'e}termin{\'e} modulo
$\R'$.

Le lemme r{\'e}sulte de ce que l'espace vectoriel $V_{\alpha'}$
form{\'e} des {\'e}l{\'e}ments de la forme $\sum_{j=1}^{n-1}x_je_j$
o{\`u} $x_j$ sont des polyn{\^o}mes de degr{\'e} strictement inf{\'e}rieur {\`a} $\alpha_j$
est suppl{\'e}mentaire {\`a} tout $\R'\in S_\alpha'({\bar k})$ dans
$\bigoplus_{j=1}^{n-1}e_j{\bar{\cal O}}$. $\square$

\begin{corollaire}
$S_\alpha$ est isomorphe {\`a} l'espace affine de dimension
$$\<\alpha,(n-1,\ldots,1,0)\>.$$
\end{corollaire}

\begin{lemme}
\begin{enumerate}
\item Soient $\mu$ et $\lambda$ deux $n$-partitions de $d$ avec $\mu<\lambda$.
On a $S_\lambda\cap X_\mu=\O$.
\item L'intersection $S_\lambda\cap X_\lambda$ est un espace affine de
dimension $2\<\lambda,\delta\>$.
\item La restriction de $h_\lambda$ {\`a} $S_\lambda\cap X_\lambda$
est constante de valeur $0$.
\end{enumerate}
\end{lemme}

\DEMONSTRATION
\begin{enumerate}
\item
Soit $\R=x{\bar{\cal O}}^n\in (S_\lambda\cap X_\mu)({\bar k})$ o{\`u} $x$ 
est une matrice comme dans le lemme pr{\'e}c{\'e}dent et o{\`u} $\mu$ et $\lambda$
sont deux $n$-partitions de $d$.
Tous les mineurs d'ordre $i$ de $x$ sont alors divisibles par
$\vp^{\mu_{n-i+1}+\cdots+\mu_n}$.
En consid{\'e}rant la sous-matrice form{\'e}e des $i$ derni{\`e}res lignes et
des $i$ derni{\`e}re colonnes, on obtient l'in{\'e}galit{\'e}
$$\lambda_{n-i+1}+\cdots+\lambda_n\geq \mu_{n-i+1}+\cdots+\mu_n$$
d'o{\`u} $\mu\geq\lambda$.

\item
Supposons maintenant que $\mu=\lambda$. Consid{\'e}rons la sous-matrice
$(i+1)\times(i+1)$ de $x$ incluant le coefficient $x_{j,n-i}$
avec $j<n-i$ et incluant les $i$ derni{\`e}res lignes ainsi que les $i$ derni{\`e}res
colonnes de $x$. Il r{\'e}sulte de la condition port{\'e}e sur les mineurs
que le polyn{\^o}me $x_{j,n-i}$ est divisible par $\vp^{\lambda_{n-i}}$.

Si les coefficients $x_{j,k}$ sont divisibles par
$\vp^{\lambda_k}$ pour tout $j<k$, alors
$x\in N({\bar{\cal O}})\vp^\lambda$
si bien que $x{\bar{\cal O}}^n\in(S_\lambda\cap X_\lambda)({\bar k})$.

Il s'ensuit que $S_\lambda\cap X_\lambda$ est isomorphe {\`a} l'espace
affine de dimension
$$\displaylines{
\<\lambda,(n-1,\ldots,1,0)\>-\<\lambda,(0,1,\ldots,n-1)\>\cr
        =\<\lambda,(n-1,n-3,\ldots,1-n)\>.\cr}$$

\item
On a d{\'e}montr{\'e} que
$$(S_\lambda\cap X_\lambda)({\bar k})
=N({\bar\OO})\vp^\lambda{\bar\OO}^n$$
si bien que la restriction de $h_\lambda$ {\`a} $S_\lambda\cap X_\lambda$
est constante et de valeur nulle. $\square$
\end{enumerate}

\begin{corollaire}
On a un isomorphisme
$$\RR\Gamma_c(S_\lambda\otimes_k{\bar k},
\A_\lambda\otimes h_\lambda^*\L_\psi)
=\Ql[-2\<\lambda,\delta\>](-\<\lambda,\delta\>).$$
\end{corollaire}

\DEMONSTRATION
On sait d'apr{\`e}s le lemme pr{\'e}c{\'e}dent que
$$S_\lambda\cap{\bar X}_\lambda=S_\lambda\cap X_\lambda$$
si bien que la restriction de $\A_\lambda$ {\`a} $S_\lambda$
est isomorphe {\`a} $\Ql[2\<\lambda,\delta\>](\<\lambda,\delta\>)$
support{\'e} par l'espace affine $S_\lambda\cap X_\lambda$
de dimension $2\<\lambda,\delta\>$. $\square$

\section{Etude de ${\tilde S}_\lambda$}

Posons ${\tilde S}_\lambda=S_\lambda\times_{X_d}{\tilde X}_d$.
L'ensemble des ${\bar k}$-points de ${\tilde S}_\lambda$ est
l'ensemble des drapeaux de r{\'e}seaux
$${\bar\OO}^n=\R_0\supset\R_1\supset\cdots\supset\R_d=\R$$
o{\`u} $\dim_{\bar k}(\R_{i-1}/\R_i)=1$ et o{\`u} $\R\in S_\lambda({\bar k})$.

Un tel drapeau {\'e}tant fix{\'e},
Pour chaque $i=0,\ldots,d$, il existe $\alpha^i\in\NN^n$
avec $|\alpha^i|=i$ tel que $\R_i\in S_{\alpha_i}({\bar k})$
Le sch{\'e}ma ${\tilde S}_\alpha$ est ainsi stratifi{\'e} selon la donn{\'e}e
d'une matrice
$\tau=(\alpha^i_j)^{0\leq i\leq d}_{1\leq j\leq n}\in\NN^{\,(d+1)n}$
telle que
\begin{itemize}
\item $\alpha^{i-1}_j\leq\alpha^i_j$ ;
\item $\sum_{j=1}^n\alpha^i_j=i$ ;
\item $\alpha^d=\lambda$.
\end{itemize}
Notons $S_\tau$ la strate correspondant {\`a} $\tau$. D{\'e}signons par
${\tilde h}_\lambda$ la restriction de 
$h_\lambda\circ\pi|_{{\tilde S}_\lambda}$ {\`a} $S_{\tau}$.

\begin{proposition}
S'il existe un $d'$ avec $1\leq d'\leq d-1$ tel que la suite
$(\alpha^{d'}_j)_{1\leq j\leq n}$ n'est pas d{\'e}croissante, alors on a
$$\RR\Gamma_c(S_\tau\otimes_k{\bar k},{\tilde h}_\lambda^*\L_\psi)=0.$$  
\end{proposition}

\DEMONSTRATION
Soit $\tau'=(\alpha^{i'}_j)^{0\leq i\leq d'}_{1\leq j\leq n}$
la sous-matrice form{\'e}e des $d'+1$ premi{\`e}res colonnes de $\tau$.
Notons ${\pi'}:S_\tau\rta S_{\tau'}$ le morphisme d{\'e}fini par
$${\pi'}(\R_0\supset\R_1\supset\cdots\supset\R_d)
=(\R_0\supset\R_1\supset\cdots\supset\R_{d'}).$$
On va d{\'e}montrer que $\RR{\pi'}_* {\tilde h}_\lambda^*\L_\psi=0$
ce qui implique par la suite spectrale de Leray que 
$$\RR\Gamma_c(S_\tau\otimes_k{\bar k},{\tilde h}_\lambda^*\L_\psi)=0.$$

Fixons un point g{\'e}om{\'e}trique
$$\R_\bullet=(\R_0\supset\cdots\supset \R_{d'}=\R')
\in S_{\tau'}({\bar k}).$$
Le groupe $\GL(\R')\cap N({\bar F})$
vu comme ${\bar k}$-groupe alg{\'e}brique de dimension infinie,
agit naturellement sur la fibre
$$\displaylines{\qquad {\pi'}^{-1}(\R_\bullet)=\{\R'=\R_{d'}\supset\R_{d'+1}
\supset\cdots\supset\R_d\mid \hfill \cr\hfill
\dim_{\bar k}(\R_{i-1}/\R_i)=1~{\rm et}~
\R_i\in S_{\alpha_i}({\bar k})\}.\qquad}$$

\begin{lemme}
Si $\alpha^{d'}$ n'est pas d{\'e}croissante,
pour tout $\R'\in S_{\alpha^{d'}}({\bar k})$,
il existe un sous-groupe
$${\Bbb G}_{a,{\bar k}}\subset \GL(\R')\cap N({\bar F})$$
tel que la restriction du caract{\`e}re $N({\bar F})\rta{\Bbb G}_{a,{\bar k}}$
d{\'e}fini par
$$n\mapsto\res(\sum_{i=1}^{n-1}n_{i,i+1}\d\vp)$$
{\`a} ce sous-groupe est l'identit{\'e} de ${\Bbb G}_{a,{\bar k}}$.
\end{lemme}

\noindent{\it D{\'e}monstration du lemme.}
Consid{\'e}rons d'abord le cas $\R'=\vp^{\alpha^{d'}}{\bar\OO}^n$.
Il existe un entier $j$ tel que $\alpha^{d'}_j<\alpha^{d'}_{j+1}$.
Le sous-groupe form{\'e} des {\'e}l{\'e}ments $n\in N({\bar F})$ tels que
$n_{k,l}=0$ avec $k<l$, $(k,l)\not=(j,j+1)$
et $n_{j,j+1}\in {\Bbb G}_{a,{\bar k}}\vp^{-1}\d\vp$
stabilise le r{\'e}seau $\vp^{\alpha^{d'}}{\bar\OO}^n$ et donc remplit
toutes les conditions requises par le lemme.

Si $\R=x{\bar\OO}^n$ pour un certain $x\in N({\bar F})$, il suffit
de conjuguer le ${\Bbb G}_{a,{\bar k}}$ pr{\'e}c{\'e}dent par $x$. $\carre$

\noindent{\it Fin de la d{\'e}monstration de la proposition.}
Notons $Z={\pi'}^{-1}(\R_\bullet)$ et $h$ la restriction de
${\tilde h}_\lambda$ {\`a} $Z$. On a une action $\xi$ de $\Ga$ sur $Z$ tel que
$h(\xi(t,z))=t+h(z)$ pour tout $t\in{\bar k}$ et $z\in Z({\bar k})$.
En particulier, on a un isomorphisme
$$\xi^*h^*\L_\psi{\tilde\rta}h^*\L_\psi\boxtimes\L_\psi.$$
La proposition r{\'e}sulte de lemme g{\'e}n{\'e}ral suivant qui est d{\'e}j{\`a}
implicite dans \cite{Del}.

\begin{lemme}
Soit $Z$ un sch{\'e}ma de type fini sur ${\bar k}$, muni d'une
action $\xi:\Ga\times Z\rta Z$. Soit ${\cal F}$ un complexe born{\'e}
sur $Z$ muni d'un isomorphisme $\xi^*{\cal F}=\L_\psi\boxtimes{\cal F}$.
Alors on a $\RR\Gamma_c(Z,{\cal F})=0$.
\end{lemme}

\DEMONSTRATION
Consid{\'e}rons le diagramme commutatif
$$\diagramme{
\Ga\times Z     & \hfl{\Xi}{} &\Ga\times Z \cr
\vfl{\pr_\Ga}{} &             & \vfl{}{\pr_\Ga}\cr
 \Ga            & \hfl{}{\Id} & \Ga \cr
}$$
o{\`u}
$$\Xi(t,z)=(t,\xi(t,z))$$
est un isomorphisme. L'isomorphisme
$$\Xi^*(\Ql\boxtimes{\cal F}){\tilde\rta}{\L_\psi\boxtimes{\cal F}}$$
induit par adjonction un isomorphisme
$$\Ql\boxtimes{\cal F}{\tilde\rta}\Xi_*(\L_\psi\boxtimes{\cal F})$$
et donc un isomorphisme
$$\Ql\boxtimes\RR\Gamma_c(Z,{\cal F}){\tilde\rta}
\L_\psi\boxtimes\RR\Gamma_c(Z,{\cal F})$$
lequel ne peut exister que si $\RR\Gamma_c(Z,{\cal F})=0$. $\square$

\begin{proposition}
Si pour tout $i=0,\ldots,d$, $\alpha^i$ est une suite d{\'e}croissante alors
on a un isomorphisme
$$\displaylines{\RR\Gamma_c(S_\tau\otimes_k{\bar k},{\tilde h}_\lambda^*\L_\psi)\cr
=\Ql[-2\<\lambda,(n-1,\ldots,1,0)\>](-\<\lambda,(n-1,\ldots,1,0)\>).}$$
\end{proposition}

\DEMONSTRATION
La proposition r{\'e}sulte du lemme suivant.

\begin{lemme}
\begin{enumerate}
\item
Pour tout $\tau$,
$S_\tau$ est isomorphe {\`a} un espace affine de dimension
$$\<\lambda,(n-1,\ldots,1,0)\>.$$
\item
Si de plus, pour tout $i$, $\alpha^i$ est une suite d{\'e}croissante alors
la restriction de ${\tilde h}_\lambda$ {\`a} $S_\tau$ est constante {\`a}
l'image nulle.
\end{enumerate}
\end{lemme}

\DEMONSTRATION
\begin{enumerate}
\item
Pour tout $i=1,\ldots,d$, vu les contraintes port{\'e}es sur les
$\alpha^i_j$, il existe un unique $j$ tel que
$\alpha^i_j=\alpha^{i-1}_j+1$. On peut en fait voir $\tau$ comme une
application
$$\{1,2,\ldots,d\}\rta\{1,2,\ldots,n\}$$
telle que pour tout $j=1,\ldots,n$, on a $|\tau^{-1}(j)|=\lambda_j$.

On d{\'e}montre par r{\'e}currence sur $d$ que $S_\tau$ est isomorphe {\`a} un
espace affine de dimension
$$\sum_{i=1}^{d}(n-\tau(i))=\sum_{j=1}^n\lambda_j(n-j).$$

Notons $\tau'$ la matrice
$(\alpha^i_j)^{0\leq i\leq d-1}_{1\leq j\leq n}\in\NN^{dn}$.
Supposons que $S_{\tau'}$ est isomorphe {\`a} un espace affine de
dimension
$$\sum_{i=1}^{d-1}(n-\tau(i))$$

Notons ${\cal F}$ le fibr{\'e} vectoriel de rang $n$ dont la fibre
au-dessus d'un point
$$\R'_\bullet=(\R_0\supset\R_1\supset\cdots\supset\R_{d-1}=\R')\in
S_{\tau'}({\bar k})$$
est l'espace vectoriel $\R'/\vp\R'$.

On peut {\'e}crire de mani{\`e}re unique
$\R'=x'{\bar\OO}^n$ avec une matrice triangulaire sup{\'e}rieure $x'$
v{\'e}rifiant les conditions de l'{\'e}nonc{\'e} du lemme 2.1. En particulier, on
a 
$${\cal F}_{\R'_\bullet}=\bigoplus_{i=1}^n \epsilon_i{\bar k}$$
o{\`u} $\epsilon_i$ est la r{\'e}duction de
$e_i\vp^{\alpha_i}+\sum_{j=1}^{i-1} x_{j,i}e_j$ modulo $\vp\R'$
si bien que le fibr{\'e} ${\cal F}$ est en fait un fibr{\'e} trivial.

De plus, la donn{\'e}e d'un ${\bar k}$-point $\R_\bullet$ de $S_\tau$
au-dessus de $\R'_\bullet$ est {\'e}quivalente {\`a} la donn{\'e}e d'un
sous-espace vectoriel de codimension $1$ de ${\cal F}_{\R'_\bullet}$ 
qui contient
$$\epsilon_1{\bar k}\oplus\cdots\oplus\epsilon_{\tau(d)-1}{\bar k}$$
mais qui ne contient pas
$$\epsilon_1{\bar k}\oplus\cdots\oplus\epsilon_{\tau(d)}{\bar k}.$$

Ce sous-espace vectoriel s'{\'e}crit de mani{\`e}re unique sous la forme 
$$\bigoplus_{j=1}^{\tau(d)-1}\epsilon_j{\bar k}
\oplus(x_{\tau(d)+1}\epsilon_{\tau(d)}+\epsilon_{\tau(d)+1}){\bar k}
\oplus\cdots\oplus
(x_{n}\epsilon_{\tau(d)}+\epsilon_{n}){\bar k}$$
si bien qu'on a un isomorphisme
$$S_{\tau'}\times\Ga^{n-\tau(d)}{\tilde\rta}S_{\tau}.$$

Compte tenu de l'hypoth{\`e}se de r{\'e}currence,
$S_\tau$ est isomorphe {\`a} un espace affine de dimension
$\sum_{i=1}^d(n-\tau(i))$.

\item
Supposons que toutes les suites $\alpha^{i}$
sont d{\'e}croissantes. On d{\'e}montre par r{\'e}currence sur $d$ que si 
$$\R_\bullet=({\bar\OO}^n=\R_0\supset\R_1\supset\cdots\supset\R_d=\R)
\in S_{\tau}({\bar k})$$
alors 
$$\R\in\N({\bar\OO})\vp^{\lambda}{\bar\OO}^n.$$

Par r{\'e}currence, on peut supposer que 
$$\R'=\R_{d-1}\in N({\bar\OO})\vp^{\alpha'}{\bar\OO}^{n-1}$$
o{\`u} $\alpha'=\alpha^{d-1}$
et quitte {\`a} utiliser l'action de $N({\bar\OO})$, on peut en fait
supposer que 
$$\R'=\vp^{\alpha'}{\bar\OO}^n.$$

Notons $l=\tau(d)$. On peut {\'e}crire 
$$\R=\bigoplus_{j=1}^{l}\vp^{\lambda_j}e_j{\bar\OO}\oplus
\bigoplus_{j=l+1}^{n}(\vp^{\lambda_j}e_j+x_j\vp^{\lambda_l-1}e_l){\bar\OO}$$
avec $x_j\in{\bar k}$ pour $j=l+1,\ldots,n$.
Du fait que $\lambda_l-1\geq \lambda_j$ pour tout $j=l+1,\ldots,n$,
on a  $\R\in N({\bar\OO})\vp^\lambda{\bar\OO}^n$. $\square$
\end{enumerate}

\noindent{\it Fin de la d{\'e}monstration du th{\'e}or{\`e}me 2. }
Pour terminer la d{\'e}monstration, il suffit de montrer que les matrices 
$\tau=(\alpha^i_j)^{0\leq i\leq d}_{1\leq j\leq n}\in\NN^{\,(d+1)n}$
telles que
\begin{itemize}
\item $\alpha^{i-1}_j\leq\alpha^i_j$ ;
\item $\sum_{j=1}^n\alpha^i_j=i$ ;
\item $\alpha^d=\lambda$ ;
\item $\alpha^{i}_{j-1}\geq\alpha^{i}_j$.
\end{itemize}
sont en correspondance univoque avec les $\lambda$-tableaux standards.

On a vu que les $\tau$ v{\'e}rifiant les trois premi{\`e}res
conditions et ne v{\'e}rifiant pas obligatoirement la quatri{\`e}me
peuvent {\^e}tre vus comme une application
$$\tau:\{1,\ldots,d\}\rta\{1,\ldots,n\}$$
telle que pour tout $j$, on a $|\tau^{-1}(j)|=\lambda_j$.
Etant donn{\'e}e une telle application, on peut inscrire successivement 
$1,2,\ldots,d$ dans le diagramme de Young $\lambda$ en {\'e}crivant le
nombre $i$ dans la premi{\`e}re case encore vide de la $j=\tau(i)$-{\`e}me
ligne.

Un tel tableau est standard si et seulement si 
$$\alpha^{i}_{j-1}\geq\alpha^{i}_j$$
pour tous les  $i=1,\ldots,d$ et $j=1,\ldots,n$.

On peut aussi raisonner de mani{\`e}re plus directe comme suit. 
L'espace vectoriel $V_\lambda$ admet une base index{\'e}e par l'ensemble
des composantes irr{\'e}ductibles de dimension maximale de la fibre
de $\pi:{\tilde X}_d\rta X_d$ au-dessus d'un point g{\'e}om{\'e}trique de
$X_\lambda$ par exemple de $\vp^\lambda{\bar\OO}^n\in X_\lambda({\bar k})$. 
En utilisant les lemmes 3.2 et 3.5, on voit facilement
que ces composantes sont pr{\'e}cis{\'e}ment les fibres des $S_\tau$ 
au-dessus de $\vp^\lambda{\bar\OO}^n$ pour les
$\tau$ dont toutes les suites $\alpha^{i}$ sont d{\'e}croissantes.

\bigskip
\noindent{\bf Remerciement} Je voudrais exprimer ma profonde gratitude
envers G{\'e}rard Laumon qui, par ses encouragements, m'a constamment
soutenu.

\def\refname{R{\'e}f{\'e}rences}

\bigskip
Ng{\^o} Bao Ch{\^a}u\\
INSTITUT GALIL{\'E}E\\
av. J.-B. Cl{\'e}ment\\
93430 Villetaneuse\\
FRANCE\\
ngo@math.univ-paris13.fr\\
 
\end{document}